\documentclass[11pt,leqno,a4paper]{article}

\usepackage{amsmath, amssymb, amsthm, amsfonts}
\usepackage{mathrsfs}
\usepackage{color}
\theoremstyle{plain}
\newtheorem{thm}{Theorem}

\newtheorem{lem}{Lemma}

\theoremstyle{remark}
\newtheorem{rem}{Remark}

\def\K{\mathop{\mbox{\bf\Large K}}}

\numberwithin{equation}{section}

\begin{document}

\title{Some new inequalities for the gamma function}
\author{Xiaodong Cao }
\date{}

\maketitle
\footnote[0]{Address: Department of Mathematics and Physics, Beijing Institute of Petro-Chemical Technology, Beijing, 102617, PR China. \\
\quad E-mail: caoxiaodong@bipt.edu.cn \\
\quad Tel/Fax:(+86)010-81292176}



\footnote[0]{This work is supported by
the National Natural Science Foundation of China (Grant No.11171344) and the Natural
Science Foundation of Beijing (Grant No.1112010).}

\begin{abstract} In this paper, we present some new inequalities for the gamma function. The main tools are the multiple-correction method developed in~\cite{CXY,Cao} and a generalized Mortici's lemma.
\end{abstract}
\section{Introduction}

Duo to its importance in mathematics, the problem of finding new and sharp inequalities for the gamma function and, in particular for large values of $x$
\begin{align}
\Gamma(x):=\int_{0}^{\infty}t^{x-1}e^{-t}dt,\quad x>0,
\end{align}
has attracted the attention of many researchers 
(see~\cite{Alz,Bat,CL,CL-1,HV,Mor1,Mor2,Mor3,Mor4,Mor5} and the references therein). Let's recall some of the classical results.
%Here we give a very brief overview of the classical results.
Maybe one of the most well-known formula for approximation the gamma function is the Stirling's formula
\begin{align}
  \Gamma(x+1)\sim \sqrt{2\pi x} \left(\frac xe\right)^{x},\quad x\rightarrow +\infty.
\end{align}
See, e.g. [1, p. 253]. The following two formulas give slightly better estimates than Stirling's formula,
\begin{align}
 &\Gamma(x+1)\sim \sqrt{2\pi } \left(\frac {x+\frac 12}{e}\right)^{x+\frac 12},\quad (\mbox{Burnside~\cite{Bur}, 1917 }),\\
&\Gamma(x+1)\sim \sqrt{2\pi } \left(\frac xe\right)^{x}\sqrt{x+\frac 16},\quad (\mbox{Gosper~\cite{Gos}, 1978}).
\end{align}
Ramanujian~\cite{Ram} proposed the claim (without proof) for
 the gamma function
\begin{align}
\Gamma (x+1)=\sqrt{\pi}\left(\frac x e\right)^x
\left(8x^3+4x^2+x+\frac{\theta_x}{30}\right)^{\frac 16},
\label{Ramanujan-C}
\end{align}
where $\theta_x\rightarrow 1$ as $x\rightarrow +\infty$ and
$\frac{3}{10}< \theta_x < 1$. This open problem was
solved by Karatsuba\cite{Kar}. Thus~\eqref{Ramanujan-C} provides a more accurate
estimate for the gamma function (see Sec. 2 below).

In this paper, we will continue the previous works~\cite{CXY,Cao}, and introduce a class of new approximations to improve these inequalities.

Throughout the paper, the notation $\Psi(k;x)$
denotes a polynomial of degree $k$ in $x$ with all coefficients non-negative, which may be different at each occurrence.
Let $(a_n)_{n\ge 1}$ and $(b_n)_{n\ge 0}$ be two sequences of real numbers with $a_n\neq 0 $ for all $n\in\mathbb{N}$ . The generalized continued fraction
\begin{align*}
\tau=b_0+\frac{a_1}{b_1+\frac{a_2}{b_2+\ddots}}=
b_0+\frac{a_1}{b_1+}\frac{a_2}{b_2+}
\cdots
=b_0+\K_{n=0}^{\infty}
\frac{a_n}{b_n}
\end{align*}
is defined as the limit of the $n$th approximant
\begin{align*}
\frac{A_n}{B_n}=b_0+\K_{k=1}^{n}\frac{a_k}{b_k}
\end{align*}
as $n$ tends to infinity. See [2, p.105].

\section{ A generalized Mortici's lemma}
Mortici~\cite{Mor1} established a very useful tool for measuring the rate of convergence, which says that a sequence $(x_n)_{n\ge 1}$ converging to zero is the fastest possible when the difference $(x_n-x_{n+1})_{n\ge 1}$ is the fastest possible. Since then, Mortici's lemma has been effectively applied in many paper such as~\cite{CXY,Cao,Mor4,Mor5}.
The following lemma is a generalization of Mortici's lemma.
\begin{lem} If $\lim_{x\rightarrow+\infty}f(x)=0$, and there exists the limit
\begin{align}
 \lim_{x\rightarrow+\infty}x^\lambda\left(f(x)-f(x+1)\right)=l\in
 \mathbb{R},
\end{align}
with $\lambda>1$, then there exists the limit
\begin{align}
 \lim_{x\rightarrow+\infty}x^{\lambda-1}f(x)=\frac{l}{\lambda-1}.
\end{align}
\end{lem}
\proof It is not very difficult to prove that for $x>2$
\begin{align}
\frac{1}{(\lambda-1)x^{\lambda-1}}=\int_{x}^{+\infty}\frac{dt}
{t^\lambda}
\le\sum_{j=0}^{\infty}\frac{1}{(x+j)^\lambda} \le
 \int_{x-1}^{+\infty}\frac{dt}{t^\lambda}
=\frac{1}{(\lambda-1)(x-1)^{\lambda-1}}.\label{SI-inequality}
\end{align}
For $\varepsilon>0$, we assume that $l-\varepsilon\le x^\lambda
\left(f(x)-f(x+1)\right)\le l+\varepsilon$ for every real number
$x$ greater than or equal to the rank $X_0>0$. By adding the
inequalities of the form
\begin{align}
(l-\varepsilon)\frac{1}{x^\lambda}
\le f(x)-f(x+1)\le (l+\varepsilon)\frac{1}{x^\lambda},
\end{align}
we get
\begin{align}
(l-\varepsilon)\sum_{j=0}^{m-1}\frac{1}{(x+j)^\lambda}
\le f(x)-f(x+m)\le (l+\varepsilon)\sum_{j=0}^{m-1}\frac{1}{(x+j)^\lambda}
\end{align}
for every $x\ge X_0$ and $m\ge 1$. By taking the limit as $m\rightarrow\infty$, then multiplying by $x^{\lambda-1}$, we obtain
\begin{align}
(l-\varepsilon)x^{\lambda-1}\sum_{j=0}^{\infty}
\frac{1}{(x+j)^\lambda}
\le x^{\lambda-1}f(x) \le (l+\varepsilon)x^{\lambda-1}\sum_{j=0}^{\infty}
\frac{1}{(x+j)^\lambda}. \end{align}
It follows from~\eqref{SI-inequality} that
\begin{align}
\frac{l-\varepsilon}{\lambda-1}\le x^{\lambda-1}f(x) \le \frac{l+\varepsilon}{\lambda-1}\frac{x^{\lambda-1}}
{(x-1)^{\lambda-1}}.
\end{align}
Now by taking the limit as $x\rightarrow+\infty$, this completes
the proof of the lemma at once.\qed

\noindent {\bf {An example}} Let's consider the Ramanujan's asymptotic formula (1.5). Let the error term $E(x)$ be defined by the following relation
\begin{align}
\Gamma (x+1)=\sqrt{\pi}\left(\frac x e\right)^x
\left(8x^3+4x^2+x+\frac{1}{30}\right)^{\frac 16}\left(1+E(x)\right).
\label{Ramanujan-C-1}
\end{align}
It follows readily from the recurrence formula $\Gamma(x+1)=x\Gamma(x)$ that
\begin{align}
\ln\left(1+E(x)\right)-\ln\left(1+E(x+1)\right)
=&-1+x\ln\left(1+\frac 1x\right)\\
&+\frac 16\ln\frac{8(x+1)^3+4(x+1)^2+(x+1)+\frac{1}{30}}
{8x^3+4x^2+x+\frac{1}{30}}.\nonumber
\end{align}
By using the \emph{Mathematica} software, we expand the right-hand function in the above formula as a power series in terms of $1/x$:
\begin{align}
\ln\left(1+E(x)\right)-\ln\left(1+E(x+1)\right)=
\frac{11}{2880 x^5}+O(\frac{1}{x^6}).
\end{align}
Thus, by Lemma 1 we have
\begin{align}
\lim_{x\rightarrow+\infty}x^4 \ln\left(1+E(x)\right)=
\frac{11}{11520}.
\end{align}
Noting that $\lim_{u\rightarrow 0}\frac{\ln(1+u)}{u}=1$, one get finally
\begin{align}
\lim_{x\rightarrow+\infty}x^4 E(x)=
\frac{11}{11520}.
\end{align}

\begin{rem}
Just as Motici's lemma, Lemma 1 also provides a method for finding the limit of a function as $x$ tends to infinity.
\end{rem}

\section{Gosper-type inequalities}
In this section, we use an example to illustrate the idea of this paper. To this end, we introduce some class of correction function $(\mathrm{MC}_k(x))_{k\ge 0}$ such that the relative error function $E_k(x)$ has the fastest possible rate of convergence, which are defined by the relations
\begin{align}
\Gamma(x+1)=\sqrt{2\pi}\left(\frac xe\right)^x \sqrt{x+\frac 16+\mathrm{MC}_k(x)}\cdot\exp(E_k(x)).
\end{align}
If $\lim_{x \rightarrow +\infty}x^{\mu}f(x)=l\neq 0$ with constant $\mu> 0$, we say that the function $f(x)$ is order $x^{-\mu}$, and write the exponent of convergence $\mu=\mu(f(x))$. Clearly if
$\mu(E_k(x))=\mu_k$, we have the following asymptotic formula
\begin{align}
\Gamma(x+1)=\sqrt{2\pi}\left(\frac xe\right)^x \sqrt{x+\frac 16+\mathrm{MC}_k(x)}\cdot
 \left(1+O(x^{-\mu_k})\right),\quad x\rightarrow +\infty.
 \end{align}
Let us briefly review a so-called \emph{multiple-correction method} presented in our previous paper~\cite{CXY,Cao}.
Actually, the \emph{multiple-correction method} is a recursive algorithm, and one of its advantages is that by repeating correction process we always can accelerate the convergence, i.e. the sequence $(\mu(E_k(x)))_{k\ge 0}$ is a strictly increasing.
The key step is to find a suitable structure of $\mathrm{MC}_k(x)$. In general, the correction function $\mathrm{MC}_k(x)$ is a finite
\emph{generalized} continued fraction (see~\cite{Cao} or~\eqref{MC-1-Def} below) or a \emph{hyper-power} series (see~\cite{CXY} or~\eqref{Ramanujan-2-2} below) in $x$.

It is not difficult to see that (3.1) is equivalent to
\begin{align}
\ln\Gamma(x+1)=\frac{1}{2}\ln (2\pi)+x\left(\ln x-1\right)+\frac 12\ln\left(x+\mathrm{MC}_k(x)\right)+E_k(x).
\end{align}
By the recurrence formula $\Gamma (x+1)=x\Gamma(x)$, we have for $x>0$
\begin{align}
E_k(x)-E_k(x+1)=-1+x\ln\left(1+\frac 1x\right)+\frac 12\ln\frac{(x+1)+\frac 16+\mathrm{MC}_k(x+1)}{x+\frac 16+\mathrm{MC}_k(x)}.
\end{align}
Now by taking the initial-correction function $\mathrm{MC}_0(x)=\frac{\kappa_0}{x+\lambda_0}$ and using \emph{Mathematica} software, we expand $E_k(x)-E_k(x+1)$ into a power series in terms of $1/x$:
\begin{align}
E_0(x)-E_0(x+1)=&\frac{-\frac{1}{72}+\kappa_0}{x^3}+\frac{17- 945\kappa_0-810 \kappa_0 \lambda_0}{540 x^4}+\\
&\frac{-641+ 33120\kappa_0-12960\kappa_0^2+43200\kappa_0 \lambda_0+25920 \kappa_0 \lambda_0^2}{12960x^5}+O\left(\frac{1}{x^6}\right).\nonumber
\end{align}
The fastest possible function $E_0(x)-E_0(x+1)$ is obtained when the first two coefficients in the above formula vanish. In this case, we find $\kappa_0=\frac {1}{72}, \lambda_0=\frac{31}{90}$, and
\begin{align}
E_0(x)-E_0(x+1)=\frac{5929}{1166400 x^5}+O\left(\frac{1}{x^6}\right).
\end{align}
By Lemma 1, we can check that
\begin{align}
\lim_{x\rightarrow +\infty}x^4E_0(x)=\frac{5929}{4665600}.
\end{align}
We continue the above correction process to successively determine the correction function $\mathrm{MC}_k(x)$ until some $k^*$ you want. On one hand, to find the related coefficients, we often use an appropriate symbolic computations software because it's huge of computations. On the other hand, the exact expressions at each occurrence also need lot of space. Hence in this paper we omit many related details. For interesting readers, see our previous paper~\cite{CXY,Cao}. In fact, we can prove that for $0\le k \le 3$
\begin{align}
\mathrm{MC}_k(x)=\K_{j=0}^{k}\frac{\kappa_j}{x+\lambda_j},
\label{MC-1-Def}
\end{align}
where
\begin{align*}
&\kappa_0=\frac{1}{72},& \lambda_0=\frac{31}{90},\\
&\kappa_1=\frac{5929}{32400},& \lambda_1=\frac{481937}{3735270},\\
&\kappa_2=\frac{76899172249}{248039857296},& \lambda_2=\frac{7745462509019287}{19149278075101482},\\
&\kappa_3=\frac{786873417270631211749921}{851541507731717527392144},
& \lambda_3=\frac{2098335745817751685364201067279071}
{30311088872486921466334781589254970}.
\end{align*}
By Lemma 1 again, we get for some constant $C_k\neq 0$
\begin{align}
\lim_{x\rightarrow +\infty}x^{2k+4}E_k(x)=C_k, \quad (k=0,1,2,3),
\end{align}
i.e. $\mu(E_k(x))=2k+4$ for $k=0,1,2,3$. Thus we obtain more accurate approximation formulas: 
\begin{align}
\Gamma(x+1)=\sqrt{2\pi}\left(\frac xe\right)^x \sqrt{x+\frac 16+\mathrm{MC}_k(x)}\cdot
 \left(1+O(x^{-(2k+4)})\right),\quad x\rightarrow +\infty.
 \end{align} 

It should be noted that if we rewrite $\mathrm{MC}_k(x)$ in the form of $\frac{P_r(m)}{Q_s(m)}$, where $P, Q$ are polynomials with $r=k$ and $s=k+1$, theoretically at least, for a large $x$ the above formula may reduce or eliminate numerically computations compared with the previous results, see e.g.~\cite{CL-1, HV}. This is the main advantage of the \emph{multiple-correction method}.

The following theorem tells us how to obtain sharp inequalities.
\begin{thm} Let $\mathrm{MC}_k(x)$ be defined as~\eqref{MC-1-Def}.
 Let $x\ge 1$, then we have for $k=0,2,$%$n\in\mathbb{N}$
\begin{align}
\Gamma(x+1)>\sqrt{2\pi}\left(\frac xe\right)^x \sqrt{x+\frac 16+\mathrm{MC}_k(x)}, \label{Gos1-1}
\end{align}
and for $k=1,3,$
\begin{align}
\Gamma(x+1)<\sqrt{2\pi}\left(\frac xe\right)^x \sqrt{x+\frac 16+\mathrm{MC}_k(x)}.\label{Gos1-2}
\end{align}
\end{thm}
\proof
We let $f_k(x)=E_k(x)-E_k(x+1)$. Clearly if $\lim_{x\rightarrow+\infty} E_k(x)=0$, then $E_k(x)=\sum_{j=0}^{\infty}f_k(x+j)$. This transformation plays an important role in this paper (essentially, it is a difference method). Hence, in order to prove inequality $E_k(x)>0 ~(\mbox {or}~E_k(x)<0)$, it suffices to show that the equality
$f_k(x)>0 ~(\mbox {or}~f_k(x)<0)$ holds under the condition $\lim_{x\rightarrow+\infty} E_k(x)=0$. By the Stirling's formula (1.2), we can show that the condition $\lim_{x\rightarrow+\infty} E_k(x)=0$ always holds. In what follows, we will apply this condition many times.

By using \emph{Mathematica} software, we may prove that for $x\ge 1$
\begin{align*}
&f_0''(x)=\frac{\Psi_1(8;x)}{x(1+x)^2(31+90x)^2(121+90x)^2(77+552x +1080 x^2)^2(1709+2712x+1080x^2)^2}>0,\\
&f_1''(x)=-\frac{\Psi_2(13;x)(x-1)+1463\cdots 9447}{x(1+x)^2)(1359251+2829648x+5976432x^2)^2\Psi_3(16;x) }<0,\\
&f_2''(x)=\frac{\Psi_4(20;x)}{x(1+x)^2\Psi_5(28;x)}>0 ,\\
&f_3''(x)=\frac{\Psi_6(25;x)(x-1)+17135\cdots 66999 }{x(1+x)^2\Psi_7(36;x)}<0.
\end{align*}
We only give the proof of inequalities in case $k=3$, other may be
proved similarly. In this case, we see that for $x\ge 1$ the inequality \eqref{Gos1-2} is equivalent to $E_3(x)<0$. As $\lim_{x\rightarrow+\infty} E_3(x)=0$, it suffices to prove that $f_3(x)<0$ for $x\ge 1$.
Since $f_3'(x)$ is strictly decreasing, but $\lim_{x\rightarrow+\infty} f_3'(x)=0$, so $f_3'(x)>0$. Thus  $f_3(x)$ is strictly increasing with $\lim_{x\rightarrow+\infty} f_3(x)=0$, so $f_3(x)<0$. This completes the proof of Theorem 1.\qed

By the \emph{multiple-correction} method, we also find another kind of inequalities.
\begin{thm} Let the $k$-th correction function $\mathrm{MC}_k(x)$ be defined by
\begin{align*}
&\mathrm{MC}_0(x)=\frac{\kappa_0}{(x+\frac{23}{90})^2+\lambda_0},\\
&\mathrm{MC}_k(x)=\frac{\kappa_0}{(x+\frac{23}{90})^2+\lambda_0+
}\K_{j=1}^{k}\frac{\kappa_j}{x+\lambda_j},\quad (k\ge 1),
\end{align*}
where
\begin{align*}
&\kappa_0=-\frac{1}{144}, &\lambda_0=\frac{4007}{21600},
\\
&\kappa_1 =\frac{4394}{637875},
&\lambda_1=\frac{130311599}{15575040},\\
&\kappa_2=\frac{7894414898425}
{119793516544},
&\lambda_2=-\frac{265702682899837009577}{34427631789478287360},\\
&\kappa_3=\frac{1897560849252106177858465792}
{77174813342532578267347147395},&
\lambda_3 =\frac{30320380455616293004898928163131563244811979}
{6134364315672065325746652708240298034227200}.
\end{align*}
Then we have
\begin{align}
 \Gamma(x+1)<\sqrt{2\pi } \left(\frac xe\right)^{x}
 \sqrt{x+\frac 16}\left(1+\mathrm{MC}_0(x)\right),\quad x\ge 13,
\end{align}
\begin{align}
 \Gamma(x+1)<\sqrt{2\pi } \left(\frac xe\right)^{x}
 \sqrt{x+\frac 16}\left(1+\mathrm{MC}_2(x)\right),\quad x\ge 6,
\end{align}
and for $k=1,3,$
\begin{align}
 \Gamma(x+1)>\sqrt{2\pi } \left(\frac xe\right)^{x}
 \sqrt{x+\frac 16}\left(1+\mathrm{MC}_k(x)\right),\quad x\ge 1.
\end{align}
\end{thm}
\proof Since the proof of Theorem 2 is very similar to that of Theorem 1, here we only give the outline of the proof. First, let the relative error function $E_k(x)$ be defined by
\begin{align}
 \Gamma(x+1)=\sqrt{2\pi } \left(\frac xe\right)^{x}
 \sqrt{x+\frac 16}\left(1+\mathrm{MC}_k(x)\right)\exp(E_k(x)).
\end{align}
Hence
\begin{align}
E_k(x)-E_k(x+1)=-1+x\ln\left(1+\frac 1x\right)+\ln\frac{1+ \mathrm{MC}_k(x+1)}{1+\mathrm{MC}_k(x)}.
\end{align}
By making use of \emph{Mathematica} software and Lemma 1, we can prove
\begin{align}
\mu(E_k(x))=2k+5,\quad (k=0,1,2,3).
\end{align}
Next we let $g_k(x)=E_k(x)-E_k(x+1)$. By using \emph{Mathematica} software, it isn't difficult to check that
\begin{align*}
&g_0''(x)=\frac{\Psi_1(14;x)(x-13)+29707\cdots 81369}{x(1+x)^2(1+6x)^2(7+6x)^2\Psi_2(16;x)}>0,\quad x\ge 13,\\
&g_1''(x)=-\frac{\Psi_3(20;x)(x-1)+13798\cdots 89479}{x(1+x)^2(1+6x)^2(7+6x)^2\Psi_4(24;x)}<0,\quad x\ge 1,\\
&g_2''(x)=\frac{\Psi_5(26;x)(x-6)+97250\cdots 34321}{x(1+x)^2(1+6x)^2(7+6x)^2\Psi_6(32;x)}>0,\quad x\ge 6,\\
&g_3''(x)=-\frac{\Psi_7(32;x)(x-1)+836559\cdots 37479}{x(1+x)^2(1+6x)^2(7+6x)^2\Psi_8(40;x)}<0,\quad x\ge 1.
\end{align*}
Lastly, just as the proof of Theorem1, Theorem 2 follows from the above inequalities readily.\qed

\section{ Ramanujan-type inequalities}
\begin{thm} Let the $k$-th correction function $\mathrm{MC}_k(x)$
be defined as
\begin{align}
\mathrm{MC}_k(x)=\K_{j=0}^{k}\frac{a_j}{x+b_j},\label{Ramanujan-2-1}
\end{align}
where
\begin{align*}
&a_0=-\frac{11}{240},&b_0=\frac{79}{154},\\
&a_1=\frac{459733}{711480},&b_1=-\frac{1455925}{70798882},
\\
&a_2=\frac{49600874140433}{101450127018720},&b_2=\frac
{10259108965771635091}{19545564575317443762},\\
&a_3=\frac{169085305336152527131511003963}{
101221579151797375403194730976},
&b_3=-\frac{6141448535908002711219920016488834171}{
203275987838924050801436670299517447102}.
\end{align*}
Let $x\ge 1$, then for $k=0,2,$
\begin{align}
\Gamma (x+1)<\sqrt{\pi}\left(\frac x e\right)^x
\left(8x^3+4x^2+x+\frac{1}{30}+\mathrm{MC}_k(x)\right)^{\frac 16},
\end{align}
and for $k=1,3,$
\begin{align}
\Gamma (x+1)>\sqrt{\pi}\left(\frac x e\right)^x
\left(8x^3+4x^2+x+\frac{1}{30}+\mathrm{MC}_k(x)\right)^{\frac 16}.
\end{align}
\end{thm}
\proof We define the relative error function $E_k(x)$ by the relation
\begin{align}
\Gamma (x+1)=\sqrt{\pi}\left(\frac x e\right)^x
\left(8x^3+4x^2+x+\frac{1}{30}+\mathrm{MC}_k(x)\right)
^{\frac 16}\exp(E_k(x)).
\end{align}
Thus
\begin{align}
E_k(x)-E_k(x+1)=&-1+x\ln\left(1+\frac 1x\right)\\
&+\frac 16\ln\frac{8(x+1)^3+4(x+1)^2+(x+1)+\frac{1}{30}+\mathrm{MC}_k(x+1)}
{8x^3+4x^2+x+\frac{1}{30}+\mathrm{MC}_k(x)}.\nonumber
\end{align}
By using \emph{Mathematica} software and Lemma 1, we can check
\begin{align}
\mu(E_k(x))=2k+6,\quad (k=0,1,2,3).
\end{align}
We let $U_k(x)=E_k(x)-E_k(x+1)$. By making use of \emph{Mathematica} software again, we can prove
\begin{align*}
&U_0''(x)=\frac{\Psi_1(13;x)(x-1)+
416838558509297754261614731715717}{3x(1+x)^2
(79+154x)^2(233+154x)^2\left(\Psi_{21}(3;x)(x-1)+
363565\right)^2
\Psi_{22}^2(4;x)}<0,\\
&U_1''(x)=\frac{\Psi_3(19;x)(x-1)+85653\cdots 25001}{x(1+x)^2\Psi_4(28;x)}>0,\\
&U_2''(x)=-\frac{\Psi_5(25;x)(x-1)+32968\cdots 13479}{x(1+x)^2\Psi_6(36;x)}<0,\\
&U_3''(x)=\frac{\Psi_7(31;x)(x-1)+17145\cdots 57723}{3x(1+x)^2\Psi_8(44;x)}>0.
\end{align*}
Similar to the proof of Theorem 1, we can get the desired assertions from the above inequalities.\qed

\begin{thm} Let the first-correction function $\mathrm{MC}_1^*(x)$
be defined by
\begin{align}
\mathrm{MC}_1^*(x)=\frac{\kappa_0}{x+\lambda_0}+
\frac{\kappa_1}{x^3+\lambda_{10}x^2+\lambda_{11}x+\lambda_{12}},
\label{Ramanujan-2-2}
\end{align}
where
\begin{align*}
&\kappa_0=-\frac{11}{240},&\lambda_0=\frac{79}{154},\\
&\kappa_1=\frac{459733}{15523200}, &\lambda_{10}=\frac{71181889}{70798882},\\
&\lambda_{11}=\frac{717183502490887}{520777318696096}, &\lambda_{12}=\frac
{1118629052995381153799}{1958878792277282473920}.
\end{align*}
Then for $x\ge 1$, the following inequality holds true
\begin{align}
\Gamma (x+1)<\sqrt{\pi}\left(\frac x e\right)^x
\left(8x^3+4x^2+x+\frac{1}{30}+\mathrm{MC}_1^*(x)\right)^{\frac 16}.\label{Ramanujan-2-2-bound}
\end{align}
\end{thm}
\proof
Let the first-correction error function $E_1^*(x)$ be defined by
\begin{align}
\Gamma (x+1)=\sqrt{\pi}\left(\frac x e\right)^x
\left(8x^3+4x^2+x+\frac{1}{30}+\mathrm{MC}_1^*(x)\right)
^{\frac 16}\exp(E_1^*(x)).
\end{align}
Hence
\begin{align}
E_1^*(x)-E_1^*(x+1)=&-1+x\ln\left(1+\frac 1x\right)\\
&+\frac 16\ln\frac{8(x+1)^3+4(x+1)^2+(x+1)+\frac{1}{30}+\mathrm{MC}^*_1(x+1)}
{8x^3+4x^2+x+\frac{1}{30}+\mathrm{MC}_1^*(x)}.\nonumber
\end{align}
By using \emph{Mathematica} software and Lemma 1, we have
\begin{align}
\mu(E_1^*(x))=10.
\end{align}
Now we let $V(x)=E_1^*(x)-E_1^*(x+1)$. By using \emph{Mathematica} again, we have
\begin{align}
V_1''(x)=-\frac{\Psi_1(33;x)(x-1)+96057\cdots 27429}{3x\left(\Psi_2(3;x)\right)^2\Psi_3(12;x)\left(\Psi_4(6;x)(x-1)
+2169\cdots 3461\right)^2 \Psi_5(14;x)}<0.\label{VD2}
\end{align}
By the same approach as the proof of Theorem 1, the inequality \eqref{Ramanujan-2-2-bound} follows from the \eqref{VD2}.\qed

\begin{rem} It is an interesting question whether our method may be used to obtain some sharp bounds for the ratio of the gamma functions, see e.g.~\cite{GXQ,Qi,Qi-1,QL}.
\end{rem}

\begin{flushleft}
Xiaodong Cao\\
Department of Mathematics and Physics, \\
Beijing Institute of Petro-Chemical Technology,\\
Beijing, 102617, P. R. China \\
E-mail: caoxiaodong@bipt.edu.cn
\end{flushleft}
\end{document}